# A GENERALIZATION OF A LEIBNIZ GEOMETRICAL THEOREM

Mihály Bencze, Florin Popovici,
Department of Mathematics, Áprily Lajos College, Braşov, Romania

Florentin Smarandache, Chair of Department of Math & Sciences, University of New Mexico, 200 College Road, Gallup, NM 87301, USA, E-mail: smarand@unm.edu

**Abstract**:
In this article we present a generalization of a Leibniz's theorem in geometry and an application of this.

**Leibniz's theorem.** Let $M$ be an arbitrary point in the plane of the triangle $ABC$, then $MA^2 + MB^2 + MC^2 = \frac{1}{3}(a^2 + b^2 + c^2) + 3MG^2$, where $G$ is the centroid of the triangle. We generalize this theorem:

**Theorem.** Let's consider $A_1, A_2, \ldots, A_n$ arbitrary points in space and $G$ the centroid of this points system; then for an arbitrary point $M$ of the space is valid the following equation:

$$\sum_{i=1}^{n} MA_i^2 = \frac{1}{n} \sum_{1 \le i < j \le n} A_i A_j^2 + n \cdot MG^2.$$

**Proof.** First, we interpret the centroid of the $n$ points system in a recurrent way.

If $n = 2$ then is the midpoint of the segment.

If $n = 3$, then it is the centroid of the triangle.

Suppose that we found the centroid of the $n-1$ points created system. Now we join each of the $n$ points with the centroid of the $n-1$ points created system; and we obtain $n$ bisectors of the sides. It is easy to show that these $n$ medians are concurrent segments. In this manner we obtain the centroid of the $n$ points created system. We'll denote $G_i$ the centroid of the $A_k$, $k = 1, 2, \ldots, i-1, i+1, \ldots, n$ points created system. It can be shown that $(n-1)A_iG = GG_i$. Now by induction we prove the theorem.

If $n = 2$ the $MA_1^2 + MA_2^2 = \frac{1}{2}A_1A_2^2 + 2MG^2$

or

$$MG^2 = \frac{1}{4}\left(2\left(MA_1^2 + MA_2^2\right)\right),$$

where $G$ is the midpoint of the segment $A_1A_2$. The above formula is the side bisector's formula in the triangle $MA_1A_2$. The proof can be done by Stewart's theorem, cosines theorem, generalized theorem of Pythagoras, or can be done vectorial. Suppose that the assertion of the theorem is true for $n = k$. If $A_1, A_2, ..., A_k$ are arbitrary points in space, $G_0$ is the centroid of this points system, then we have the following relation:

$$\sum_{i=1}^{k} MA_i^2 = \frac{1}{k} \sum_{1 \leq i < j \leq k} A_i A_j^2 + k \cdot MG_0^k.$$

Now we prove for $n = k+1$.

Let $A_{k+1} \notin \{A_1, A_2, ..., A_k, G_0\}$ be an arbitrary point in the space and let $G$ be the centroid of the $A_1, A_2, ..., A_k, A_{k+1}$ points system. Taking into account that $G$ is on the segment $A_{k+1}G_0$ and $k \cdot A_{k+1}G = GG_0$, we apply Stewart's theorem to the points $M$, $G_0$, $G$, $A_{k+1}$, from where:

$$MA_{k+1}^2 \cdot GG_0 + MG_0^2 \cdot GA_{k+1} - MG^2 \cdot A_{k+1}G_0 = GG_0 \cdot GA_{k+1} \cdot A_{k+1}G_0.$$

According to the previous observation $A_{k+1}G = \frac{k}{k+1} A_{k+1}G_0$

and $GG_0 = \frac{k}{k+1} A_{k+1}G_0$.

Using these, the above relation becomes:

$$MA_{k+1}^2 + k \cdot MG_0^2 = \frac{k}{k+1} A_{k+1}G_0^2 + (k+1)MG^2.$$

From here

$$k \cdot MG_0^2 = \sum_{i=1}^{k} MA_i^2 - \frac{1}{k} \sum_{1 \leq i < j \leq k} A_i A_j^2.$$

From the supposition of the induction, with $M \equiv A_{k+1}$ as substitution, we obtain

$$\sum_{i=1}^{k} A_i A_j^2 = \frac{1}{k} \sum_{1 \leq i < j \leq k} A_i A_j^2 + k \cdot A_{k+1}G_0^2$$

and thus

$$\frac{k}{k+1} A_{k+1}G_0^2 = \frac{1}{k+1} \sum_{i=1}^{k} A_i A_{k+1}^2 - \frac{1}{k(k+1)} \sum_{1 \leq i < j \leq k} A_i A_j^2.$$

Substituting this in the above relation we obtain that

$$\sum_{i=1}^{k+1} MA_i^2 = \left(\frac{1}{k} - \frac{1}{k(k+1)}\right) \sum_{1 \leq i < j \leq k} A_i A_j^2 + \frac{1}{k+1} \sum_{i=1}^{k} A_i A_{k+1}^2 + (k+1)MG^2 =$$

$$= \frac{1}{k+1} \sum_{1 \leq i < j \leq k+1} A_i A_j^2 + (k+1)MG^2.$$

With this we proved that our assertion is true for $n = k+1$. According to the induction, it is true for every $n \geq 2$ natural numbers.

**Application 1.** If the points $A_1, A_2, ..., A_n$ are on the sphere with the center $O$ and radius $R$, then using in the theorem the substitution $M \equiv O$ we obtain the identity:

$$OG^2 = R^2 - \frac{1}{n^2} \sum_{1 \leq i < j \leq n} A_i A_j^2.$$

In case of a triangle: $OG^2 = R^2 - \frac{1}{9}(a^2 + b^2 + c^2)$.

In case of a tetrahedron: $OG^2 = R^2 - \frac{1}{16}(a^2 + b^2 + c^2 + d^2 + e^2 + f^2)$.

**Application 2.** If the points $A_1, A_2, ..., A_n$ are on the sphere with the center $O$ and radius $R$, then $\sum_{1 \leq i < j \leq n} A_i A_j^2 \leq n^2 R^2$.

The equality holds if and only if $G \equiv O$. In case of a triangle: $a^2 + b^2 + c^2 \leq 9R^2$, in case of a tetrahedron: $a^2 + b^2 + c^2 + d^2 + e^2 + f^2 \leq 16R^2$.

**Application 3.** Using the arithmetic and harmonic mean inequality, from the previous application, it results the following inequality:

$$\sum_{1 \leq i < j \leq n} \frac{1}{A_i A_j^2} \geq \frac{(n-1)^2}{4R^2}.$$

In the case of a triangle: $\frac{1}{a^2} + \frac{1}{b^2} + \frac{1}{c^2} \geq \frac{1}{R^2}$, in case of a tetrahedron:

$$\frac{1}{a^2} + \frac{1}{b^2} + \frac{1}{c^2} + \frac{1}{d^2} + \frac{1}{e^2} + \frac{1}{f^2} \geq \frac{9}{4R^2}.$$

**Application 4.** Considering the Cauchy-Buniakowski-Schwarz inequality from the Application 2, we obtain the following inequality:

$$\sum_{1 \leq i < j \leq n} A_i A_j^2 \leq nR\sqrt{\frac{n(n-1)}{2}}.$$

In case of a triangle: $a + b + c \leq 3\sqrt{3}R$, in case of a tetrahedron:
$$a + b + c + d + e + f \leq 4\sqrt{6}R.$$

**Application 5.** Using the arithmetic and harmonic mean inequality, from the previous application we obtain the following inequality

$$\sum_{1 \leq i < j \leq n} \frac{1}{A_i A_j^2} \geq \frac{(n-1)\sqrt{n(n-1)}}{2R\sqrt{2}}.$$

In case of a triangle: $\frac{1}{a} + \frac{1}{b} + \frac{1}{c} \geq \frac{\sqrt{3}}{R}$, in case of a tetrahedron:

$$\frac{1}{a} + \frac{1}{b} + \frac{1}{c} + \frac{1}{d} + \frac{1}{e} + \frac{1}{f} \geq \frac{3}{R}\sqrt{\frac{3}{2}}.$$

**Application 6.** Considering application 3, we obtain the following inequality:

$$\frac{n^2(n-1)^2}{4} \leq \left( \sum_{1 \leq i < j \leq n} A_i A_j^k \right) \left( \sum_{1 \leq i < j \leq n} \frac{1}{A_i A_j^k} \right) \leq$$

$$\leq \begin{cases} \dfrac{(M+m)^2 n^2(n-1)^2}{16M \cdot m} & \text{if } \dfrac{n(n-1)}{2} \text{ is even,} \\ \dfrac{(M+m)^2 n^2(n-1)^2 - 4(M-m)^2}{16M \cdot m} & \text{if } \dfrac{n(n-1)}{2} \text{ is odd} \end{cases}$$

where $m = \min \{A_i A_j^k\}$ and $M = \max \{A_i A_j^k\}$. In case of a triangle:

$$9 \leq (a^k + b^k + c^k)(a^{-k} + b^{-k} + c^{-k}) \leq \frac{2M^2 + 5M \cdot m + 2m^2}{M \cdot m},$$

in case of a tetrahedron:

$$36 \leq (a^k + b^k + c^k + d^k + e^k + f^k)(a^{-k} + b^{-k} + c^{-k} + d^{-k} + e^{-k} + f^{-k}) \leq \frac{9(M+m)^2}{M \cdot m}.$$

**Application 7.** Let $A_1, A_2, \ldots, A_n$ be the vertexes of the polygon inscribed in the sphere with the center $O$ and radius $R$. First we interpret the orthocenter of the inscribable polygon $A_1 A_2 \ldots A_n$. For three arbitrary vertexes, corresponds one orthocenter. Now we take four vertexes. In the obtained four orthocenters of the triangles we construct the circles with radius $R$, which have one common point. This will be the orthocenter of the inscribable quadrilateral. We continue in the same way. The circles with radius $R$ that we construct in the orthocenters of the $n-1$ sides inscribable polygons have one common point. This will be the orthocenter of the $n$ sides, inscribable polygon. It can be shown that $O$, $H$, $G$ are collinear and $n \cdot OG = OH$. From the first application

$$OH^2 = n^2 R^2 - \sum_{1 \leq i < j \leq n} A_i A_j^2$$

and

$$GH^2 = (n-1)^2 R^2 - \left(1 - \frac{1}{n}\right)^2 \sum_{1 \leq i < j \leq n} A_i A_j^2.$$

In case of a triangle $OH^2 = 9R^2 - (a^2 + b^2 + c^2)$ and $GH^2 = 4R^2 - \dfrac{4}{9}(a^2 + b^2 + c^2)$.

**Application 8.** In the case of an $A_1 A_2 \ldots A_n$ inscribable polygon $\sum_{1 \leq i < j \leq n} A_i A_j^2 = n^2 R^2$ if and only if $O \equiv H \equiv G$. In case of a triangle this is equivalent with an equilateral triangle.

**Application 9.** Now we compute the length of the midpoints created by the $A_1, A_2, \ldots, A_n$ space points system. Let $S = \{1, 2, \ldots, i-1, i+1, \ldots, n\}$ and $G_0$ be the centroid of the $A_k$, $k \in S$, points system. By substituting $M \equiv A_i$ in the theorem, for the length of the midpoints we obtain the following relation:

$$A_i G_0^2 = \frac{1}{n-1} \sum_{k \in S} A_i A_k^2 - \frac{1}{(n-1)^2} \sum_{u, v \in S : u \neq v} A_u A_v^2.$$

**Application 10.** In case of a triangle $m_a^2 = \dfrac{2(b^2+c^2)-a^2}{4}$ and its permutations.

From here:
$$m_a^2 + m_b^2 + m_c^2 = \frac{3}{4}(a^2+b^2+c^2),$$
$$m_a^2 + m_b^2 + m_c^2 \le \frac{27}{4}R^2,$$
$$m_a + m_b + m_c \le \frac{9}{2}R.$$

**Application 11.** In case of a tetrahedron $m_a^2 = \dfrac{1}{9}\left(3(a^2+b^2+c^2)-(d^2+e^2+f^2)\right)$ and its permutations.

From here:
$$\sum m_a^2 = \frac{4}{9}\left(\sum a^2\right),$$
$$\sum m_a^2 \le \frac{64}{9}R^2,$$
$$\sum m_a \le \frac{16}{3}R.$$

**Application 12.** Denote $m_{a,f}$ the length of the segments, which join midpoint of the $a$ and $f$ skew sides of the tetrahedron (bimedian). In the interpretation of the application $9m_{a,f}^2 = \dfrac{1}{4}\left(b^2+c^2+d^2+e^2-a^2-f^2\right)$ and its permutations.

From here
$$m_{a,f}^2 + m_{b,d}^2 + m_{c,e}^2 = \frac{1}{4}\left(\sum a^2\right),$$
$$m_{a,f}^2 + m_{b,d}^2 + m_{c,e}^2 \le 4R^2,$$
$$m_{a,f} + m_{b,d} + m_{c,e} \le 2R\sqrt{3}.$$